\documentclass[12pt]{amsart}
\usepackage{amssymb}
\usepackage{amsfonts}
\usepackage{latexsym}
\usepackage{amscd}
\usepackage[mathscr]{euscript}
\usepackage{xy} \xyoption{all}
\numberwithin{equation}{section}

\newcommand{\Z}{{\mathbb{Z}}}
\newcommand{\N}{{\mathbb{N}}}

\newcommand{\C}{{\mathbb{C}}}
\newcommand{\uloopr}[1]{\ar@'{@+{[0,0]+(-4,5)}@+{[0,0]+(0,10)}@+{[0,0] +(4,5)}}^{#1}}
\newcommand{\uloopd}[1]{\ar@'{@+{[0,0]+(5,4)}@+{[0,0]+(10,0)}@+{[0,0]+ (5,-4)}}^{#1}}
\newcommand{\dloopr}[1]{\ar@'{@+{[0,0]+(-4,-5)}@+{[0,0]+(0,-10)}@+{[0, 0]+(4,-5)}}_{#1}}
\newcommand{\dloopd}[1]{\ar@'{@+{[0,0]+(-5,4)}@+{[0,0]+(-10,0)}@+{[0,0 ]+(-5,-4)}}_{#1}}

\newcommand{\luloop}[1]{\ar@'{@+{[0,0]+(-8,2)}@+{[0,0]+(-10,10)}@+{[0, 0]+(2,2)}}^{#1}}

\newtheorem{lem}{Lemma}[section]
\newtheorem{corol}[lem]{Corollary}
\newtheorem{theor}[lem]{Theorem}
\newtheorem{prop}[lem]{Proposition}
\newtheorem{rema}[lem]{Remark}
\newtheorem{defi}[lem]{Definition}

\newtheorem{exems}[lem]{Examples}

\begin{document}
\title[Stable rank of Leavitt path algebras]{Stable rank of Leavitt path algebras}
\author{P. Ara}
\address{Departament de Matem\`atiques, Universitat Aut\`onoma de Barcelona,
08193 Bellaterra (Barcelona), Spain.} \email{para@mat.uab.cat}
\author{E. Pardo}
\address{Departamento de Matem\'aticas, Universidad de C\'adiz,
Apartado 40, 11510 Puerto Real (C\'adiz),
Spain.}\email{enrique.pardo@uca.es}
\urladdr{http://www2.uca.es/dept/matematicas/PPersonales/PardoEspino/index.HTML}
\thanks{The first author was partially supported by the DGI
and European Regional Development Fund, jointly, through Project
MTM2005-00934. The second author was partially supported by the DGI
and European Regional Development Fund, jointly, through Project
MTM2004-00149, and by PAI III grants FQM-298 and P06-FQM-1889 of the
Junta de Andaluc\'{\i}a. Both authors are partially supported by
by the Comissionat per Universitats i Recerca de la Generalitat de
Catalunya.} 
%
%
\begin{abstract}
We characterize the values of the stable rank for Leavitt path
algebras, by giving concrete criteria in terms of properties of the
underlying graph.
\end{abstract}

\maketitle

\section*{Introduction and background}

Leavitt path algebras have been recently introduced in \cite{AA1}
and \cite{AMFP}. Given an arbitrary (but fixed) field $K$ and a
row-finite graph $E$, the Leavitt path algebra $L_K(E)$ is the
algebraic analogue of the Cuntz-Krieger algebra $C^*(E)$ described
in \cite{Raeburn}. Several interesting ring-theoretic properties
have been characterized for this class of algebras. For instance,
the Leavitt path algebras which are purely infinite simple have been
characterized in \cite{AA2}, and \cite{APS} contains a
characterization of the Leavitt path algebras which are exchange
rings in terms of condition (K), a purely graph-theoretic condition
defined below.

In this paper, we show that the only possible values of the (Bass)
stable rank of a Leavitt path algebra are $1$, $2$ and $\infty$.
Moreover a precise characterization in terms of properties of the
graph of the value of the stable rank is provided (Theorem
\ref{Th:Gordo}). A similar result was obtained in \cite[Theorem
7.6]{APS} under the additional hypothesis that the graph satisfies
condition (K) (equivalently, $L(E)$ is an exchange ring). Many tools
of the proof of that result must be re-worked in our general
situation. We have obtained in several situations simpler arguments
that work without the additional hypothesis of condition (K) on the
graph. Another new feature of our approach is a detailed analysis of
the stable rank in extensions of Leavitt path algebras of stable
rank $2$, in order to show that the stable rank of these extensions
cannot shift to $3$. Our main tool for this study is the well-known
concept of elementary rank; see for example \cite[Chapter 11]{McR}.

The (topological) stable rank of the Cuntz-Krieger algebras $C^*(E)$
was computed in \cite{DHSz}. This paper has been the inspiration for
the general strategy of the proof here. Note that, by a result of
Herman and Vaserstein \cite{HV}, the topological and the Bass stable
ranks coincide for $C^*$-algebras. For the sake of comparison, let
us mention that, although the possible values of the stable rank of
$C^*(E)$ are also $1$, $2$ and $\infty$, it turns out that there are
graphs $E$ such that the stable rank of $C^*(E)$ is $1$ while the
stable rank of $L(E)$ is $2$. Concretely, if $E$ is a graph such
that no cycle has an exit and $E$ contains some cycle, then the
stable rank of $C^*(E)$ is $1$ by \cite[Theorem 3.4]{DHSz}, but the
stable rank of $L(E)$ is $2$ by Theorem \ref{Th:Gordo}.

Along this paper, we describe the Leavitt path algebras following
the presentation of \cite[Sections 3 and 5]{AMFP}, but using the
notation of \cite{AA1} for the elements.

A \emph{(directed) graph} $E=(E^0,E^1,r,s)$ consists of two
countable sets  $E^0,E^1$ and maps $r,s:E^1 \to E^0$. The elements
of $E^0$ are called \emph{vertices} and the elements of $E^1$
\emph{edges}.

A vertex which emits no edges is called a  \emph{sink}. A graph $E$
is \emph{finite} if $E^0$ and $E^1$ are finite sets.  If $s^{-1}(v)$
is a  finite set for every $v\in E^0$, then the graph is called
\emph{row-finite}. A \emph{path} $\mu$ in a graph $E$ is a sequence
of edges $\mu=(\mu_1, \dots, \mu_n)$ such that
$r(\mu_i)=s(\mu_{i+1})$ for $i=1,\dots,n-1$. In such a case,
$s(\mu):=s(\mu_1)$ is the \emph{source} of $\mu$ and
$r(\mu):=r(\mu_n)$ is the \emph{range} of $\mu$. If $s(\mu)=r(\mu)$
and $s(\mu_i)\neq s(\mu_j)$ for every $i\neq j$, then $\mu$  is a
called a \emph{cycle}. We say that a cycle $\mu =(\mu_1, \dots
,\mu_n) $ has an \emph{exit} if there is a vertex $v=s(\mu_i)$ and
an edge $f\in s^{-1}(v)\setminus \{\mu_i\}$. If $v=s(\mu)=r(\mu)$
and $s(\mu_i)\neq v$ for every $i>1$, then $\mu$ is a called a
\emph{closed simple path based at $v$}. We denote by $CSP_E(v)$ the
set of closed simple paths in $E$ based at $v$. For a path $\mu$ we
denote by $\mu^0$ the set of its vertices, i.e.,
$\{s(\mu_1),r(\mu_i)\mid i=1,\dots,n\}$. For $n\ge 2$ we define
$E^n$ to be the set of paths of length $n$, and $E^*=\bigcup_{n\ge
0} E^n$ the set of all paths. We define a relation $\ge$ on $E^0$ by
setting $v\ge w$ if there is a path $\mu\in E^*$ with $s(\mu)=v$ and
$r(\mu)=w$. A subset $H$ of $E^0$ is called \emph{hereditary} if
$v\ge w$ and $v\in H$ imply $w\in H$. A set is \emph{saturated} if
every vertex which feeds into $H$ and only into $H$ is again in $H$,
that is, if $s^{-1}(v)\neq \emptyset$ and $r(s^{-1}(v))\subseteq H$
imply $v\in H$. Denote by $\mathcal{H}$ (or by $\mathcal{H}_E$ when
it is necessary to emphasize the dependence on $E$) the set of
hereditary saturated subsets of $E^0$. We denote by $E^\infty$ the
set of infinite paths $\gamma=(\gamma_n)_{n=1}^\infty$ of the graph
$E$ and by $E^{\le \infty}$ the set $E^\infty$ together with the set
of finite paths in $E$ whose end vertex is a sink. We say that a
vertex $v$ in a graph $E$ is \emph{cofinal} if for every $\gamma\in
E^{\le \infty}$ there is a vertex $w$ in the path $\gamma$ such that
$v\ge w$. We say that a  graph $E$ is \emph{cofinal} if so are all
the vertices of $E$. According to \cite[Lemma 2.8]{APS}, this is
equivalent to the fact that $\mathcal{H}=\{ \emptyset , E^0 \}$.

Let $E=(E^0,E^1, r, s)$ be a  graph, and let $K$ be a field. We
define the {\em Leavitt path algebra} $L_K(E)$ associated with $E$
as the $K$-algebra generated by a set $\{v\mid v\in E^0\}$ of
pairwise orthogonal idempotents, together with a set of variables
$\{e,e^*\mid e\in E^1\}$, which satisfy the following relations:

(1) $s(e)e=er(e)=e$ for all $e\in E^1$.

(2) $r(e)e^*=e^*s(e)=e^*$ for all $e\in E^1$.

(3) $e^*e'=\delta _{e,e'}r(e)$ for all $e,e'\in E^1$.

(4) $v=\sum _{\{ e\in E^1\mid s(e)=v \}}ee^*$ for every $v\in E^0$
that emits edges.

Note that the relations above imply that $\{ee^*\mid e\in E^1\}$ is
a set of pairwise orthogonal idempotents in $L_K(E)$. Note also that
if $E$ is a finite graph then we have $\sum _{v\in E^0} v=1$. In
general the algebra $L_K(E)$ is not unital, but it can be written as
a direct limit of unital Leavitt path algebras (with non-unital
transition maps), so that it is an algebra with local units (recall
that a \emph{local unit} in a ring $R$ is an increasing net of
idempotents $\{e_{\lambda}\}_{\lambda \in \Lambda}\subset R$ such
that for each $a\in R$ there exists $\mu \in \Lambda$ with
$a=ae_{\mu}=e_{\mu}a$). Along this paper, we will be concerned only
with row-finite graphs $E$ and we will work with Leavitt path
algebras over an arbitrary but fixed field $K$. We will usually
suppress the field from the notation.

Recall that $L(E)$ has a $\mathbb{Z}$-grading. For every $e\in E^1$,
set the degree of $e$ as $1$, the degree of $e^\ast$ as $-1$, and
the degree of every element in $E^0$ as 0. Then we obtain a
well-defined degree on the Leavitt path $K$-algebra $L(E)$, thus,
$L(E)$ is a $\mathbb{Z}$-graded algebra:
$$L(E)=\bigoplus\limits_{n\in \mathbb{Z}}L(E)_n, \quad L(E)_nL(E)_m\subseteq
L(E)_{n+m}, \  \hbox{for all }\ n, m \in \mathbb{Z}.$$

An ideal $I$ of a $\mathbb{Z}$-graded algebra $A=\oplus_{n\in
\mathbb{Z}}A_n$ is a \emph{graded ideal} in case $I=\oplus_{n\in
\mathbb{Z}} (I\cap A_n)$. By \cite[Proposition 5.2 and Theorem
5.3]{AMFP}, an ideal $J$ of $L(E)$ is graded if and only if it is
generated by idempotents; in fact, $J$ is a graded ideal if and only
if $J$ coincides with the ideal $I(H)$ of $L(E)$ generated by $H$,
where $H=J\cap E^0\in \mathcal{H}_E$. Indeed, the map $H\mapsto
I(H)$ defines a lattice isomorphism between $\mathcal{H}_E$ and
$\mathcal{L}_{\text{gr}}(L(E))$.

Recall that a graph $E$ satisfies condition (L) if every closed
simple path has an exit, and satisfies condition (K) if for each
vertex $v$ on a closed simple path there exists at least two
distinct closed simple paths $\alpha, \beta$ based at $v$.

Section 1 contains some basic information on the structure of
Leavitt path algebras, which will be very useful for the
computations in Section 2 of the stable rank of such algebras.
Finally, Section 3 contains some illustrative examples of Leavitt
path algebras.

\section{Basic facts}

For a graph $E$ and a hereditary subset $H$ of $E^0$, we denote by
$E_H$ the {\it restriction graph}
$$(H, \{e\in E^1\mid s(e)\in H\}, r|_{(E_H)^1},s|_{(E_H)^1}).$$
Observe that if $H$ is finite then $L(E_H)=p_HL(E)p_H$, where
$p_H=\sum _{v\in H} v\in L(E)$. On the other hand, for $X\in
\mathcal{H}_E$, we denote by $E/X$ the {\it quotient graph}
$$(E^0\setminus X,\{ e\in E^1\mid r(e)\notin X\}, r|_{(E/X)^1},s|_{(E/X)^1})$$
By \cite[Lemma 2.3(i)]{APS} we have a natural isomorphism
$L(E)/I(X)\cong L(E/X)$ for $X\in \mathcal H_E$. Our next result
shows that $I(X)$ is also a Leavitt path algebra.

\begin{defi}\label{assocideal}
{\rm (\cite[Definition 1.3]{DHSz}) Let $E$ be a graph, and let
$\emptyset \ne X\in \mathcal{H}_E$. Define $$F_E(X)=\{ \alpha
=(\alpha_1 \dots \alpha_n)\in E^n\mid n\ge 1,  s(\alpha _1)\in
E^0\setminus X, r(\alpha _i)\in E^0\setminus X \mbox{ for every }
i<n, r(\alpha _n)\in X\}.$$ Let
$\overline{F_E(X)}=\{\overline{\alpha}\mid \alpha \in F_E(X)\}$.
Then, we define the graph ${}_XE=({}_XE^0, {}_XE^1, s', r')$ as
follows:
\begin{enumerate}
\item ${}_XE^0=X\cup F_E(X)$. \item ${}_XE^1=\{ e\in E^1\mid
s(e)\in X\}\cup \overline{F_E(X)}$. \item For every $e\in E^1 \mbox{
with } s(e)\in X$, $s'(e)=s(e)$ and $r'(e)=r(e)$. \item For every
$\overline{\alpha}\in \overline{F_E(X)}$,
$s'(\overline{\alpha})=\alpha$ and
$r'(\overline{\alpha})=r(\alpha)$.
\end{enumerate}}
\end{defi}

\begin{lem}\label{isoideal}
Let $E$ be a graph, and let $\emptyset \ne X\in \mathcal{H}_E$.
Then, $I(X)\cong L({}_XE)$ (as nonunital rings).
\end{lem}
\begin{proof}
We define a map $\phi: L({}_XE)\rightarrow I(X)$ by the following
rule: (i) For every $v\in X$, $\phi (v)=v$; (ii) For every $\alpha
\in F_E(X)$, $\phi (\alpha )=\alpha \alpha ^*$; (iii) For every
$e\in E^1 \mbox{ with } s(e)\in X$, $\phi (e)=e$ and $\phi
(e^*)=e^*$; (iv) For every $\overline{\alpha}\in \overline{F_E(X)}$,
$\phi (\overline{\alpha})=\alpha$ and $\phi
(\overline{\alpha}^*)=\alpha^*$.

By definition, it is clear that the images of the generators of
$L({}_XE)$ satisfy the relations defining $L({}_XE)$. Thus, $\phi$
is a well-defined $K$-algebra morphism.

To see that $\phi$ is onto, it is enough to show that every vertex
of $X$ and every finite path $\alpha$ of $E$ which ranges in $X$ are
in the image of $\phi$. For any $v\in X$, $\phi (v)=v$, so that this
case is clear. Now, let $\alpha =(\alpha_1 \dots \alpha_n)$ with
$\alpha _i\in E^1$. If $s(\alpha _1)\in X$, then $\alpha=\phi
(\alpha _{1})\cdots \phi (\alpha _n)$. Suppose that $s(\alpha _1)\in
E^0\setminus X$ and $r(\alpha _n)\in X$. Then, there exists $1\leq
j\leq n-1$ such that $r(\alpha _j)\in E^0\setminus X$ and $r(\alpha
_{j+1})\in X$. Thus, $\alpha =(\alpha_1, \dots
\alpha_{j+1})(\alpha_{j+2}, \dots \alpha_n)$, where
$\beta=(\alpha_1, \dots \alpha_{j+1})\in F_E(X)$. Hence, $\alpha=
\phi (\overline{\beta})\phi (\alpha _{j+2})\cdots \phi (\alpha _n)$.

To show injectivity, notice that, for every $\alpha \in F_E(X)$,
$\alpha =\overline{\alpha}\overline{\alpha}^*$. Hence, every element
$t\in L({}_XE)$ can be written as

\begin{equation*}
\tag{1}\label{equ:1}\qquad t=\sum\limits_{\alpha, \beta \in F_E(X)}
{\overline{\alpha}a_{\alpha, \beta}\overline{\beta}^*},
\end{equation*}
where $a_{\alpha, \beta}\in L(E_X)$. Suppose that $0\ne
\mbox{Ker}(\phi)$, and let $0\ne t\in \mbox{Ker}(\phi)$ written as
in (\ref{equ:1}). By definition of the map $\phi$,
\begin{equation*}
\tag{2}\label{equ:2} \qquad 0=\phi (t)=\sum\limits_{\alpha, \beta
\in F_E(X)}{\alpha a_{\alpha, \beta}\beta^*}.
\end{equation*}
Let $\alpha _0\in F_E(X)$ with maximal length among those appearing
(with a nonzero coefficient) in the expression (\ref{equ:2}). Then,
for any other $\alpha \in F_E(X)$ appearing in the same expression,
$\alpha _0^*\cdot\alpha$ is $0$ if $\alpha \ne \alpha_0$ or
$r(\alpha_0)$ if $\alpha =\alpha _0$. Thus,
\begin{equation*}
\tag{3}\label{equ:3} \qquad 0=\sum\limits_{\alpha, \beta \in
F_E(X)}{\alpha _0^*\alpha a_{\alpha,
\beta}\beta^*}=\sum\limits_{\beta \in F_E(X)}{a_{\alpha_0,
\beta}\beta^*}.
\end{equation*}
Now, let $\beta _0\in F_E(X)$ with maximal length among those
appearing in the expression (\ref{equ:3}). The same argument as
above shows that
\begin{equation*}
\tag{4}\label{equ:4} \qquad 0=\sum\limits_{\beta \in
F_E(X)}{a_{\alpha_0, \beta}\beta^*}\beta _0 =a_{\alpha_0, \beta_0}.
\end{equation*}
But $0\ne a_{\alpha_0, \beta_0}$ by hypothesis, and we reach a
contradiction. Thus, we conclude that $\phi$ is injective.
\end{proof}

\begin{lem}\label{quasinoquot}
Let $R$ be a ring, and let $I\lhd R$ an ideal with local unit. If
there exists an ideal $J\lhd I$ such that $I/J$ is a unital simple
ring, then there exists an ideal $M\lhd R$ such that $R/M\cong I/J$.
\end{lem}
\begin{proof}
Given $a\in J$, there exists $x\in I$ such that $a=ax=xa$. Thus,
$J\subseteq JI$, and $J\subseteq IJ$. Hence, $J\lhd R$.

By hypothesis, there exists an element $e\in I$ such that
$\overline{e}\in I/J$ is the unit. Consider the set $\mathcal{C}$ of
ideals $L$ of $R$ such that $J\subseteq L$ and $e\not\in L$. If we
order $\mathcal{C}$ by inclusion, it is easy to see that it is
inductive. Thus, by Zorn's Lemma, there exists a maximal element of
$\mathcal{C}$, say $M$. Then, $J\subseteq M\cap I\subsetneqq I$,
whence $J=M\cap I$ by the maximality of $J$ in $I$. Thus,
$$I/J=I/(M\cap I)\cong I+M/M\lhd R/M.$$ Suppose that $R\ne I+M$.
Clearly, $\overline{e}\in I+M/M$ is a unit. Thus, $\overline{e}$ is
a central idempotent of $R/M$ generating $I+M/M$. So,
$L=\{a-a\overline{e}\mid a\in R/M\}$ is an ideal of $R/M$, while
$$R/M=\overline{e}(R/M)+L,$$
being the sum an internal direct sum. If $\pi \colon
R\twoheadrightarrow R/M$ is the natural projection map, then $\pi
^{-1}(L)=M+\{a-ae\mid a\in R\}$ is an ideal of $R$ containing $M$
(and so $J$). If $e\in \pi ^{-1}(L)$, then $L=R/M$, which is
impossible. Hence, $\pi^{-1}(L)\in \mathcal{C}$, and contains
strictly $M$, contradicting the maximality of $M$ in $\mathcal{C}$.
Thus, $I+M=R$, and so $R/M\cong I/J$, as desired.
\end{proof}

\begin{corol}\label{noquot}
Let $E$ be a graph, and let $H\in\mathcal{H}_E$. If there exists
$J\lhd I(H)$ such that $I(H)/J$ is a unital simple ring, then there
exists an ideal $M\lhd L(E)$ such that $L(E)/M\cong I(H)/J$.
\end{corol}
\begin{proof}
By Lemma \ref{isoideal}, $I(H)\cong L({}_HE)$, whence $I(H)$ has a
local unit. Thus, the result holds by Lemma \ref{quasinoquot}.
\end{proof}

Recall that an idempotent $e$ in a ring $R$ is called {\it infinite}
if $eR$ is isomorphic as a right $R$-module to a proper direct
summand of itself. A simple ring $R$ is called {\it purely infinite}
in case every nonzero right ideal of $R$ contains an infinite
idempotent. See \cite{AGP} for some basic properties of purely
infinite simple rings and \cite[Theorem 11]{AA2} for a
characterization of purely infinite simple Leavitt path algebras in
terms of properties of the graph.

\begin{prop}\label{P:pisu->gradedmax}
Let $E$ be a row-finite graph, and let $J$ be a maximal two-sided
ideal of $L(E)$. If $L(E)/J$ is a unital purely infinite simple
ring, then $J\in \mathcal{L}_{\text{gr}}(L(E))$.
\end{prop}
\begin{proof} Let $a$ be an element of $L(E)$ such that $a+J$ is
the unit in $L(E)/J$. There are $v_1,\cdots ,v_n\in E^0$ such that
$a\in pL(E)p$, where $p=v_1+\cdots +v_n\in L(E)$. Since $av=va=0$
for all $v\in E^0\setminus \{v_1,\dots ,v_n\}$, it follows that the
hereditary saturated set $X=\{ v\in E^0\mid v\in J\}$ is cofinite in
$E^0$ and thus passing to $L(E)/I(X)\cong L(E/X)$, we can assume
that $E$ is a finite graph and that $E^0\cap J=\emptyset$.

Since $E$ is finite, the lattice $\mathcal{L}_{\text{gr}}(L(E))$ of
graded ideals (equivalently, idempotent-generated ideals) of $L(E)$
is finite by \cite[Theorem 5.3]{AMFP}, so that there exists a
nonempty $H\in \mathcal{H}_E$ such that $I=I(H)$ is minimal as a
graded ideal. Since $I+J=L(E)$ by our assumption that $J\cap
E^0=\emptyset$, we have
$$I/(I\cap J)\cong L(E)/J,$$
so that $I$ has a unital purely infinite simple quotient. Since
$I\cong L(_HE)$  and $J\cap I$ does not contain nonzero idempotents,
it follows from our previous argument that $_HE$ is finite and so
$I$ is unital. So $I=eL(E)$ for a central idempotent $e$ in $L(E)$.
Since $I$ is graded-simple, \cite[Remark 6.7]{AMFP} and
\cite[Theorem 11]{AA2} imply that $I$ is either $M_n(K)$ or
$M_n(K[x,x^{-1}])$ for some $n\ge 1$ or it is simple purely
infinite. Since $I$ has a quotient algebra which is simple purely
infinite, it follows that $I\cap J=0$ and $J=(1-e)L(E)$ is a graded
ideal. Indeed we get $e=1$ because we are assuming that $J$ does not
contain nonzero idempotents.
\end{proof}

Notice that, as a consequence of Proposition \ref{P:pisu->gradedmax}
and \cite[Theorem 11]{AA2}, we get the following generalization of
\cite[Lemma 7.2]{APS}, which is analog to \cite[Proposition
3.1]{DHSz}

\begin{lem}\label{L:pisuquot}
Let $E$ be a row-finite graph. Then, $L(E)$ has a unital purely
infinite simple quotient if and only if there exists $H\in
\mathcal{H}_E$ such that the quotient graph $E/H$ is nonempty,
finite, cofinal, contains no sinks and each cycle has an exit.
\end{lem}

\section{Stable rank for Leavitt path algebras}

Let $S$ be any unital ring containing an associative ring $R$ as a
two-sided ideal. The following definitions can be found in
\cite{Vas}. A column vector $b=(b_i)_{i=1}^n$ is called
\emph{$R$-unimodular} if $b_1-1,b_i\in R$ for $i>1$ and there exist
$a_1-1,a_i\in R$ ($i>1$) such that $\sum_{i=1}^n a_ib_i=1$. The
\emph{stable rank} of $R$ (denoted by $\mbox{sr}(R)$) is the least
natural number $m$ for which for any $R$-unimodular vector
$b=(b_i)_{i=1}^{m+1}$ there exist $v_i\in R$ such that the vector
$(b_i+v_ib_{m+1})_{i=1}^m$ is $R$-unimodular. If such a natural $m$
does not exist we say that the stable rank of $R$ is infinite.

\begin{lem}\label{rankone} {\rm (cf. \cite[Lemma 7.1]{APS})}
Let $E$ be an acyclic graph. Then, the stable rank of $L(E)$ is $1$.
\end{lem}

\begin{lem}\label{rankinf}
Let $E$ be a graph. If there exists a unital purely infinite simple
quotient of $L(E)$, then the stable rank of $L(E)$ is $\infty$.
\end{lem}
\begin{proof}
If there exists a maximal ideal $M\lhd L(E)$ such that $L(E)/M$ is a
unital purely infinite simple ring, then $\mbox{sr}(L(E)/M)=\infty$
(see e.g. \cite{AGP}). Since $\mbox{sr}(L(E)/M)\leq \mbox{sr}(L(E))$
(see \cite[Theorem 4]{Vas}), we conclude that
$\mbox{sr}(L(E))=\infty$.
\end{proof}

We adapt the following terminology from \cite{DHSz}: we say that a
graph $E$ has isolated cycles if whenever $(a_1, \dots ,a_n)$ and
$(b_1, \dots ,b_m)$ are closed simple paths such that
$s(a_i)=s(b_j)$ for some $i,j$, then $a_i=b_j$. Notice that, in
particular, if $E$ has isolated cycles, the only closed simple paths
it can contain are cycles.

\begin{lem}\label{ranktwo} {\rm (cf. \cite[Lemma 3.2]{DHSz}, \cite[Lemma 7.4]{APS})}
Let $E$ be a graph. If $L(E)$ does not have any unital purely
infinite simple quotient, then there exists a graded ideal $J\lhd
L(E)$ with ${\mbox{sr}}(J)\le 2$ such that $L(E)/J$ is isomorphic to
the Leavitt path algebra of a graph with isolated cycles. Moreover
$\mbox{sr}(J)=1$ if and only if $J=0$.
\end{lem}
\begin{proof} Set $$X_0=\{v\in E^0\mid \exists e\ne f\in E^1 \text{ with }
s(e)=s(f)=v,\,\, r(e)\ge v,\,\, r(f)\ge v\},$$ and let $X$ be the
hereditary saturated closure of $X_0$. Consider $J=I(X)$. Then $J$
is a graded ideal of $L(E)$ and $L(E)/J\cong L(E/X)$ by \cite[Lemma
2.3(1)]{APS}. It is clear from the definition of $X_0$ that $E/X$ is
a graph with isolated cycles.

It remains to show that $\mbox{sr}(J)\le 2$ and that
$\mbox{sr}(J)=2$ if $J\ne 0$. The proof of these facts follows the
lines of the proof of \cite[Lemma 7.4]{APS}, using Corollary
\ref{noquot} instead of \cite[Proposition 5.4]{APS} and Lemma
\ref{isoideal} instead of \cite[Lemma 5.2]{APS}.
\end{proof}

\begin{defi}\label{srextens}
Let $A$ be a unital ring with stable rank $n$. We say that $A$ has
stable rank {\it closed by extensions} in case for any unital ring
extension
$$\begin{CD}
0 @>>> I @>>> B @>>> A @>>> 0
\end{CD}$$
of $A$ with $\mbox{sr}(I)\le n$ we have $\mbox{sr}(B)=n$.
\end{defi}

Recall that a unital ring $R$ has {\it elementary rank} $n$, denoted
by $\mbox{er}(R)=n$, in case that, for every $t\ge n+1$, the
elementary group $E_t(R)$ acts transitively on the set $U_c(t,R)$ of
$t$-unimodular columns with coefficients in $R$, see
\cite[11.3.9]{McR}.

In the next lemma, we collect some properties that we will need in
the sequel.

\begin{lem}\label{L:elemen}
Let $A$ be a unital ring. Assume that $\mbox{sr}(A)=n<\infty$.
\begin{enumerate}
   \item If $\mbox{er}(A)<n$ then $M_m(A)$ has stable rank closed by extensions
   for every $m\ge 1$.
   \item Let $D$ be any (commutative) euclidean domain  such that
   $\mbox{sr}(D)>1$  and let $m$ be a positive integer. Then
   $\mbox{sr}(M_m(D))=2$ and $\mbox{er}(M_m(D))=1$. In particular $M_m(D)$ has
   stable rank closed by extensions.
   \item Let $$\begin{CD} 0 @>>> I @>>> B @>>> A @>>> 0
\end{CD}$$
be a unital extension of $A$. If $\mbox{er}(A)<n$ and $I$ has a
local unit $(g_i)$ such that
   $\mbox{sr}(g_iIg_i)\le n$ and $\mbox{er}(g_iIg_i)<n$ for all $i$, then
   $\mbox{sr}(B)=n$ and $\mbox{er}(B)<n$.

\end{enumerate}
\end{lem}

\begin{proof}
(1) This is essentially contained in \cite{Vas}. We include a sketch
of the proof for the convenience of the reader. Assume that we have
a unital extension $B$ of $A$ with $\mbox{sr}(I)\le n$. Let
$\textbf{a}=(a_1,\dots,a_{n+1})^t\in U_c(n+1,B)$. Then
$\overline{\textbf{a}}=(\overline{a_1},\dots,\overline{a_{n+1}})^t\in
U_c(n+1,A)$. Since $\mbox{sr}(A)=n$, there exists $b_1,\dots,b_n\in
B$ such that $(\overline{a_1}+\overline{b_1}\overline{a_{n+1}},
\dots, \overline{a_n}+\overline{b_n}\overline{a_{n+1}})^t\in
U_c(n,A)$. Replacing $\textbf{a}$ with $(a_1+b_1a_{n+1}, \dots
,a_n+b_na_{n+1}, a_{n+1})$, we can assume that $(\overline{a_1},
\dots ,\overline{a_n} )^t\in U_c(n,A)$.

Since $\mbox{er}(A)\le n-1$, there exists $E\in E(n,B)$ such that
$\overline{E}\cdot (\overline{a_1}, \dots ,\overline{a_n} )^t=(1,
0,\dots ,0)^t$. Since $\textbf{a}$ is reducible if and only if
$\mbox{diag}(E,1)\cdot \textbf{a}$ is reducible, we can assume that
$(\overline{a_1}, \dots ,\overline{a_n} )^t=(1, 0,\dots ,0)^t$.
Finally, replacing $a_{n+1}$ with $a_{n+1}-a_1a_{n+1}$, we can
assume that $\overline{\textbf{a}}=(1, 0,\dots ,0)^t$, that is,
$\textbf{a}\in U_c(n+1, I)$. Now, as $\mbox{sr}(I)\le n$,
$\textbf{a}$ is reducible in $I$, and so in $B$, as desired.

Given any positive integer $m\ge 1$, $\mbox{sr}(M_m(A))=\lceil
(\mbox{sr}(A)-1)/m\rceil+1$ by \cite[Theorem 3]{Vas} and
$\mbox{er}(M_m(A))\le\lceil \mbox{er}(A)/m\rceil$ by \cite[Theorem
11.5.15]{McR}. So, it is clear that $\mbox{er}(A)<\mbox{sr}(A)$
implies $\mbox{er}(M_m(A))<\mbox{sr}(M_m(A))$. Hence, by the first
part of the proof, $M_m(A)$ has stable rank closed by extensions, as
desired.

(2) It is well known that an Euclidean domain has stable rank less
than or equal to $2$, and that it has elementary rank equal to $1$,
see e.g. \cite[Proposition 11.5.3]{McR}. So, the result follows from
part (1).

(3) Since $\mbox{sr}(I)\le n$, the fact that $\mbox{sr}(B)=n$
follows from part (1). Now, take $m\ge n$, and set
$\textbf{a}=(a_1,\dots, a_m)^t\in U_c(m,B)$. Since $\mbox{er}(A)<n$,
there exists $E\in E(m, B)$ such that $\overline{E}\cdot
\overline{\textbf{a}}=(1, 0, \dots ,0)^t$. So, $\textbf{b}:=E\cdot
\textbf{a}\equiv (1, 0, \dots ,0)^t (\mbox{mod }I)$. Let $g\in I$ an
idempotent in the local unit such that $b_1-1, b_2, \dots, b_m\in
gIg$. Since $\mbox{er}(gIg)< n$ by hypothesis, there exists $G\in
E(m,gIg)$ such that $(G+\mbox{diag}(1-g,\dots ,1-g))\cdot
\textbf{b}=(1, 0, \dots ,0)^t$.
\end{proof}

\begin{corol}\label{C:senseunitat}
Let $A$ be a unital $K$-algebra with $\mbox{sr}(A)=n\ge 2$ and
$\mbox{er}(A)<\mbox{sr}(A)$. Then, for any non necessarily unital
$K$-algebra $B$ and two-sided ideal $I$ of $B$ such that $B/I\cong
A$ and $\mbox{sr}(I)\le n$, we have $\mbox{sr}(B)= n$.
\end{corol}
\begin{proof}
Given any $K$-algebra $R$, we define the unitization $R^1=R\times
K$, with the product $$(r, a)\cdot (s,b)=(rs+as+rb, ab).$$ Consider
the unital extension
$$\begin{CD} 0 @>>> I @>>> B^1 @>>> A^1 @>>> 0.
\end{CD}$$
Notice that $A^1\cong A\times K$, because $A$ is unital. So,
$\mbox{sr}(A^1)=\mbox{sr}(A)$ and $\mbox{er}(A^1)=\mbox{er}(A)$.
Now, by Lemma \ref{L:elemen}(1), $\mbox{sr}(B^1)\le n$. Since $n\le
\mbox{sr}(B)\le \mbox{sr}(B^1)\le n$, the result holds.
\end{proof}

\begin{prop}\label{P:finite}
Let $E$ be a finite graph with isolated cycles. Then
$\mbox{sr}(L(E))\le 2$ and $\mbox{er}(L(E))=1$. Moreover,
$\mbox{sr}(L(E))=1$ if and only if $E$ is acyclic.
\end{prop}
\begin{proof}
We proceed by induction on the number of cycles of $E$. If $E$ has
no cycles then $\mbox{sr}(L(E))=1$ by Lemma \ref{rankone}, so that
$\mbox{er}(L(E))=1$ by \cite[Proposition 11.3.11]{McR}. Assume that
$E$ has cycles $C_1, \dots, C_n$. Define a binary relation on the
set of cycles by setting  $C_i\ge C_j$ iff there exists a finite
path $\alpha$ such that $s(\alpha)\in C_i^0$ and $r(\alpha)\in
C_j^0$. Since $E$ is a graph with isolated cycles, $\ge$ turns out
to be a partial order on the set of cycles. Since the set of cycles
is finite, there exists a maximal one, say $C_1$. Set $A=\{e\in
E^1\mid s(e)\in C_1 \mbox{ and } r(e)\not\in C_1\}$, let $S(E)$
denote the set of sinks of $E$, and define $B=\{r(e)\mid e\in
A\}\cup S(E)\cup \bigcup_{i=2}^nC_i^0$. Let $H$ be the hereditary
and saturated closure of $B$. By construction of $H$, $C_1$ is the
unique cycle in $E/H$, and it has no exits. Moreover, $E/H$
coincides with the hereditary and saturated closure of $C_1$. Then,
$L(E/H)\cong M_k(K[x, x^{-1}])$ for some $k\ge 1$. Consider the
extension
$$\begin{CD} 0 @>>> I(H) @>>> L(E) @>>> L(E/H) @>>> 0
\end{CD}.$$
Now, by Lemma \ref{L:elemen}(2), $\mbox{sr}(L(E/H))=2$ and
$\mbox{er}(L(E/H))=1$. Consider the local unit $(p_X)$ of
$L({}_HE)\cong I(H)$ consisting of idempotents $p_X=\sum_{v\in X} v$
where $X$ ranges on the set of vertices of ${}_HE$ containing $H$.
Since these sets are hereditary in $({_HE})^0$, we get that
$p_XI(H)p_X=p_XL({}_HE)p_X=L(({}_HE)_X)$ is a path algebra of a
graph with isolated cycles, containing exactly $n-1$ cycles. By
induction hypothesis, $\mbox{sr}(p_XI(H)p_X)\le 2$ and
$\mbox{er}(p_XI(H)p_X)=1$. So, by Lemma \ref{L:elemen}(3), we
conclude that $\mbox{sr}(L(E))=2$ and $\mbox{er}(L(E))=1$. Hence,
the induction step works, so we are done.
\end{proof}

We are now ready to obtain our main result.

\begin{theor}\label{Th:Gordo}
Let $E$ be a row-finite graph. Then the values of the stable rank of
$L(E)$ are:
\begin{enumerate}
\item $\mbox{sr}(L(E))=1$ if $E$ is acyclic. \item
$\mbox{sr}(L(E))=\infty$ if there exists $H\in \mathcal{H}_E$ such
that the quotient graph $E/H$ is nonempty, finite, cofinal, contains
no sinks and each cycle has an exit.
\item $\mbox{sr}(L(E))=2$ otherwise.
\end{enumerate}
\end{theor}
\begin{proof}
(1) derives from Lemma \ref{rankone}, while (2) derives from Lemma
\ref{rankinf} and Lemma \ref{L:pisuquot}. We can thus assume that
$E$ contains cycles and, using Lemma \ref{L:pisuquot}, that $L(E)$
does not have any unital purely infinite simple quotient.

By Lemma \ref{ranktwo}, there exists a hereditary saturated set $X$
of $E^0$ such that $\mbox{sr}(I(X))\le 2$, while $E/X$ is a graph
having isolated cycles. By \cite[Lemma 3.1]{AMFP}, there is an
ascending sequence $(E_i)$ of complete finite subgraphs of $E/X$
such that $E/X=\bigcup_{i\ge 1}E_i$. So, by \cite[Lemma 3.2]{AMFP},
$L(E/X)\cong \varinjlim L(E_i)$. For each $i\ge 1$, there is a
natural graded $K$-algebra homomorphism $\phi _i\colon L(E_i)\to
L(E/X)$. The kernel of $\phi _i$ is a graded ideal of $L(E_i)$ whose
intersection with $E_i^0$ is empty, so $\phi_i$ is injective and the
image $L_i$ of $L(E_i)$ through $\phi _i$ is isomorphic with
$L(E_i)$. It follows from Proposition \ref{P:finite} that, for every
$i\ge 1$, $\mbox{sr}(L_i)\le 2$ and $\mbox{er}(L_i)=1$. If $\pi:
L(E)\to L(E/X)$ denotes the natural epimorphism (see \cite[Lemma
2.3(1)]{APS}), then given any $i\ge 1$, we have
$$\begin{CD} 0 @>>> I(X) @>>> \pi^{-1}(L_i) @>>> L_i @>>> 0
\end{CD}.$$
If $\mbox{sr}(L_i)=1$, then $\mbox{sr}(\pi^{-1}(L_i))\le 2$ by
\cite[Theorem 4]{Vas}. If $\mbox{sr}(L_i)=2$ then it follows from
Corollary \ref{C:senseunitat} that $\mbox{sr}(\pi^{-1}(L_i))=2$.
Since $L(E)=\bigcup_{i\ge 1} \pi^{-1}(L_i)$ we get that
$\mbox{sr}(L(E))\le 2$. Since $E$ contains cycles we have that
either $I(X)\ne 0$ or $E/X$ contains cycles. If $I(X)\ne 0$ then
$\mbox{sr}(I(X))=2$ by Lemma \ref{ranktwo} and so
$\mbox{sr}(L(E))=2$ by \cite[Theorem 4]{Vas}. If $I(X)=0$, then $E$
has isolated cycles. Take a vertex $v$ in a cycle $C$ of $E$ and let
$H$ be the hereditary subset of $E$ generated by $v$. Observe that
$L(E_H)=pL(E)p$ for the idempotent $p=\sum_{w\in H^0} w\in \mathcal
M (L(E))$, where $\mathcal M (L(E))$ denotes the multiplier algebra
of $L(E)$; see \cite{AP}. Let $I$ be the ideal of $pL(E)p$ generated
by all the basic idempotents $r(e)$ where $e\in E^1$ is such that
$s(e)\in C$ and $r(e)\notin C$. Since $E$ has isolated cycles it
follows that $I$ is a proper ideal of $pL(E)p$ and moreover
$pL(E)p/I\cong M_k(K[x,x^{-1}])$, where $k$ is the number of
vertices in $C$. We get
$$\mbox{sr}(pL(E)p)\ge \mbox{sr}(pL(E)p/I)=2.$$
It follows that $1<\mbox{sr}(L(E))\le 2$ and thus
$\mbox{sr}(L(E))=2$, as desired.
\end{proof}

\section{Some remarks and examples}

In this section we present several examples of Leavitt path
algebras, and we compute their stable rank by using Theorem
\ref{Th:Gordo}. Several remarks on the relationship with the stable
rank of graph $C^*$-algebras are also included.

\begin{exems}\label{Ex:exemple1}
{\rm The basic examples to illustrate Theorem \ref{Th:Gordo}
coincide with those given in \cite[Example 1.4]{AA1}:
\begin{enumerate}
\item The Leavitt path algebra associated with the acyclic graph $E$
\[
{
\def\labelstyle{\displaystyle}
\xymatrix{{\bullet}^{v_1} \ar [r]  & {\bullet}^{v_2} \ar [r] &
{\bullet}^{v_3} \ar@{.}[r] & {\bullet}^{v_{n-1}} \ar [r] &
{\bullet}^{v_n}}}
\]
satisfies $L(E)\cong M_n(K)$. Thus, $\mbox{sr}(L(E))=1$ by Theorem
\ref{Th:Gordo}(1) (in this particular case, the original result is
due to Bass).\item For $n\geq 2$, the Leavitt path algebra
associated with the graph $F$
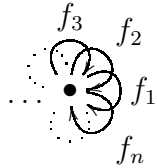
\[
{
\def\labelstyle{\displaystyle}
\xymatrix{ & {\bullet} \ar@(ur,dr) ^{f_1} \ar@(u,r) ^{f_2}
\ar@(ul,ur) ^{f_3} \ar@{.} @(l,u) \ar@{.} @(dr,dl) \ar@(r,d) ^{f_n}
\ar@{}[l] ^{\ldots} }}
\]
is an example of a unital purely infinite simple algebra, because of
\cite[Theorem 11]{AA2}; in fact $L(F)\cong L(1,n)$ --the $n$th
Leavitt algebra-- by \cite[Example 12(ii)]{AA2}. Thus,
$\mbox{sr}(L(F))=\infty $ by Theorem \ref{Th:Gordo}(2) (in this
particular case, one can also trace this fact using
\cite[Proposition 6.5]{Rieffel}).
\item Finally, the Leavitt path
algebra associated with the graph $G$
\[
{
\def\labelstyle{\displaystyle}
\xymatrix{ \bullet\dloopd{}}}
\]
satisfies $L(G)\cong K[z,z^{-1}]$ by \cite[Example 1.4(ii)]{AA1}.
Thus, $\mbox{sr}(L(G))=2 $ by Theorem \ref{Th:Gordo}(3).
\end{enumerate}
}
\end{exems}

\begin{exems}\label{Rk:notallsosimple}
{\rm We show some further examples that illustrate the complexity of
models of Leavitt path algebras:
\begin{enumerate}
\item On one hand, stable rank 2 examples can be obtained as more
or less complex extensions of the ring of Laurent polynomials, as we
can see with the Leavitt path algebra of the graph $E$
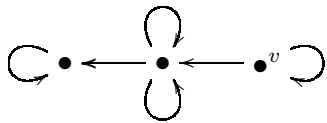
\[
{
\def\labelstyle{\displaystyle}
\xymatrix{ \bullet\dloopd{} &  \bullet \ar[l] \uloopr{}\dloopr{}
\ar[l] & \bullet^v \ar[l] \uloopd{}}}
\]
Here the ideal $I$ in Lemma \ref{ranktwo} is $I=I(E^0\setminus \{
v\})$, being $L(E)/I\cong K[x, x^{-1}]$. Notice that, because of
Lemma \ref{ranktwo}, $\mbox{sr}(I)=2$, while
$\mbox{sr}(L(E))=\mbox{sr}(L(E)/I)=2$ by Theorem \ref{Th:Gordo}(3).
The remarkable fact behind Theorem \ref{Th:Gordo} is that in the
context of Leavitt path algebras, extensions of stable rank 2 rings
by stable rank 2 ideals cannot attain stable rank 3 (in general this
is not true).
\item On the other hand, unital purely infinite simple Leavitt
path algebras turn out to be more complex than the classical Leavitt
algebras, so that there are plenty of unital Leavitt path algebras
with infinite stable rank different from the classical examples. For
example, the Leavitt path algebra of the graph $F$
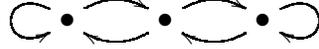
\[
{
\def\labelstyle{\displaystyle}
\xymatrix{\bullet\dloopd{} \ar@/^8pt/ [r] & {\bullet} \ar@/^8pt/ [r]
\ar@/^8pt/ [l] & {\bullet} \ar@/^8pt/ [l] \uloopd{} }}
\]
is unital purely infinite simple by \cite[Theorem 11]{AA2}, but
$(K_0(L(F)), [1])\cong (\Z, 0)$ by \cite[Theorem 3.5]{AMFP} and
\cite[Corollary 2.2]{AGP}, while $(K_0(L(1,n)), [1])\cong
(\Z/(n-1)\Z, \overline{1})$ (see \cite[Theorem 4.2]{AGP}). Taking
the graph $G$
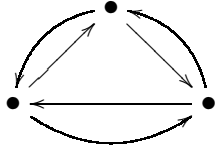
\[
{
\def\labelstyle{\displaystyle}
\xymatrix{ {} & \bullet  \ar[rd] \ar@/^{-10pt}/ [ld] &  {} \\
\bullet \ar[ru] \ar@/^{-15pt}/ [rr]&  & \bullet \ar[ll]
\ar@/^{-10pt}/ [lu] \\
}}
\]
$\mbox{ }$\vspace{.2truecm}

\noindent instead of $F$, we get a unital purely infinite simple
Leavitt path algebra such that $$(K_0(L(G)), [1])\cong
((\Z/2\Z)\oplus (\Z/2\Z), (\overline{0},\overline{0})).$$ An extra
example is that associated with the graph $H$
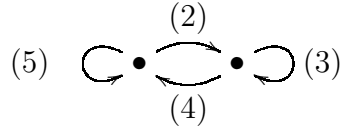
\[
{
\def\labelstyle{\displaystyle}
\xymatrix{{(5)}& \bullet\dloopd{} \ar@/^8pt/ [r] ^{(2)}& {\bullet}
\ar@/^8pt/ [l] ^{(4)}\uloopd{} ^{(3)}}}
\]
(here the $(n)$s denote the number of parallel edges), which is
again unital purely infinite simple, and such that
$$(K_0(L(H)), [1])\cong (\Z\oplus (\Z/2\Z), (1,\overline{1})).$$ No one of
them can be, then, isomorphic to any classical Leavitt algebra.
\end{enumerate}
}
\end{exems}

\begin{rema}\label{Rk:comparingSeaStar}
{\rm Fix $K=\C$ the field of complex numbers, and let $E$ be any
row-finite graph. Then:
\begin{enumerate}
\item It follows from \cite[Proposition
3.1 \& Theorem 3.4(2)]{DHSz} that $\mbox{sr}(C^*(E))=\infty$ if and
only if there exists $H\in \mathcal{H}_E$ such that the quotient
graph $E/H$ is nonempty, finite, cofinal, contains no sinks and each
cycle has an exit. By using this and Theorem \ref{Th:Gordo}, we see
that $\mbox{sr}(L(E))=\infty$ if and only if
$\mbox{sr}(C^*(E))=\infty$.
\item Since any acyclic graph is a
graph whose cycles have no exits, we have that $\mbox{sr}(L(E))=1$
implies that $\mbox{sr}(C^*(E))=1$. \item So, the only difference
occurs when $\mbox{sr}(L(E))=2$ and all the cycles in $E$ have no
exits, since then $\mbox{sr}(C^*(E))=1$ by \cite[Theorem
3.4(1)]{DHSz}. The simplest example of this situation is the graph
$G$ in Example \ref{Ex:exemple1}(3). As we noticed, $L(G)\cong \C
[z,z^{-1}]$ and $\mbox{sr}(L(G))=2$. It is clear that
$$(1+z)\C [z,z^{-1}]+ (1+z^2)\C [z,z^{-1}]=\C [z,z^{-1}],$$
and it is straightforward to see that there is no element $v\in \C
[z,z^{-1}]$ such that $(1+z)+v(1+z^2)$ is invertible in $\C
[z,z^{-1}]$. On the other hand, if we take the completion of $L(G)$,
we get the graph $C^*$-algebra $C^*(G)\cong
\mathcal{C}(\mathbb{T})$, which has stable rank 1 by
\cite[Proposition 1.7]{Rieffel}. Because of \cite{HV}, there exists
$v\in C^*(E)$ such that $(1+z)+v(1+z^2)$ is invertible in
$\mathcal{C}(\mathbb{T})$. Since a continuous function in
$\mathcal{C}(\mathbb{T})$ is invertible if and only if it has no
zeroes in $\mathbb{T}$, we see that we can take $v=1$.
\end{enumerate}
}
\end{rema}

\begin{rema}\label{Rk:Morita&Co}
{\rm Stable rank is not Morita invariant in general, but in the case
of Leavitt path algebras some interesting phenomena rise up:
\begin{enumerate}
\item Suppose that $E,F$ are finite graphs such that $L(E)$ and
$L(F)$ are Morita equivalent. Thus, $L(E)\cong P\cdot M_n(L(F))\cdot
P$ for some $n\in \N$ and some full idempotent $P\in M_n(L(F))$.
Since the values $1$ and $\infty$ in the stable rank are preserved
by passing to matrices \cite[Theorem 4]{Vas} and full corners
\cite[Theorem 7 \& Theorem 8]{AGsr}, Theorem \ref{Th:Gordo} implies
that $\mbox{sr}(L(E))=\mbox{sr}(L(F))$. So, stable rank is a Morita
invariant for unital Leavitt path algebras. \item This is not longer
true when $L(E)$ is nonunital. To see an example, let $F$ be the
graph in Example \ref{Ex:exemple1}(2), and $F^{\infty}$ be the rose
of $n$ petals
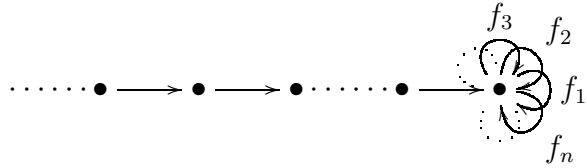
\[
{
\def\labelstyle{\displaystyle}
\xymatrix{\cdots \cdots{\bullet} \ar [r]  & {\bullet} \ar [r] &
{\bullet} \cdots \cdots {\bullet} \ar [r]  & {\bullet} \ar@(ur,dr)
^{f_1} \ar@(u,r) ^{f_2} \ar@(ul,ur) ^{f_3} \ar@{.} @(l,u) \ar@{.}
@(dr,dl) \ar@(r,d) ^{f_n}& }}
\]
with an infinite tail added. As we have seen before, $L(F)\cong
L(1,n)$ --the $n$th Leavitt algebra-- with $\mbox{sr}(L(F))=\infty$,
while an easy induction argument using \cite[Proposition 13]{AA2}
shows that $L(F^{\infty})\cong M_{\infty }(L(1,n))$. Hence these two
algebras are Morita equivalent. On the other hand, $L(F^{\infty})$
has no unital purely infinite simple quotients (as it is simple and
nonunital), so that $\mbox{sr}(L(F^{\infty}))=2$ by Theorem
\ref{Th:Gordo}(3).

Moreover, the graph $F^{\infty}$ is a direct limit (see
\cite[Section 3]{AMFP}) of the graphs $E_n^m$
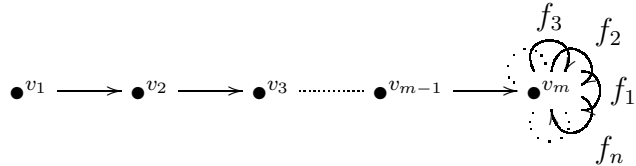
\[
{
\def\labelstyle{\displaystyle}
\xymatrix{{\bullet}^{v_1} \ar [r]  & {\bullet}^{v_2} \ar [r] &
{\bullet}^{v_3} \ar@{.}[r] & {\bullet}^{v_{m-1}} \ar [r] &
{\bullet}^{v_m} \ar@(ur,dr) ^{f_1} \ar@(u,r) ^{f_2} \ar@(ul,ur)
^{f_3} \ar@{.} @(l,u) \ar@{.} @(dr,dl) \ar@(r,d) ^{f_n}& }}
\]
introduced in \cite[Example 12]{AA2}. Since $L(F^{\infty})\cong
\varinjlim L(E_n^m)$ and $L(E_n^m)\cong M_{m}(L(1,n))$, we get
$\mbox{sr}(L(E_n^m))=\infty$ by the above remark, whence
$$2=\mbox{sr}(L(F^{\infty}))=\mbox{sr}(\varinjlim L(E_n^m))<\liminf \mbox{sr}(L(E_n^m))=\infty .$$
So, this inequality can be strict when we work with Leavitt path
algebras.
\end{enumerate}
}
\end{rema}

\section*{Acknowledgments}

Part of this work was done during visits of the first author to the
Departamento de Matem\'aticas de la Universidad de C\'adiz (Spain),
and of the second author to the Centre de Recerca Matem\`atica
(U.A.B., Spain). Both authors want to thank the host centers for
their warm hospitality. Also they wish to thank the referee for
his/her comments and suggestions.

\end{document}